\documentclass[12pt,titlepage]{amsart}
\usepackage{amsmath}

\newtheorem{thm}{Theorem}
\newtheorem{lem}[thm]{Lemma}
\newtheorem{prop}[thm]{Proposition}

\theoremstyle{definition}
\newtheorem{defn}{Definition}[section]

\theoremstyle{remark}

\theoremstyle{plain}

\newcommand{\Z}{{\mathbf Z}}

\newcommand{\R}{{\mathbf R}}

\newcommand{\supp}{\operatorname{supp}}

\newcommand{\dist}{\operatorname{dist}}

\newcommand{\nmid}{\not|}

\newcommand{\MM}{{\mathcal M}}

\newcommand{\mmod}{\operatorname{ mod }}

\newcommand{\reff}{r_{eff}}
\newcommand{\sF}{\underline F}
\newcommand{\sG}{\underline G}

\newcommand{\sO}{\underline O}
\newcommand{\sR}{\underline r}
\newcommand{\codim}{\operatorname{codim}}
\newcommand{\Af}{{\mathbf A}}
\newcommand{\CC}{{\mathcal C}}
\newcommand{\sGp}{\underline G^{(p)}}
\newcommand{\sGG}{{\mathcal G}}
\newcommand{\vol}{\operatorname{vol}}
\newcommand{\diam}{\operatorname{diam}}
\newcommand{\vl}{\vec\ell}
\newcommand{\ve}{\vec e}
\newcommand{\disc}{\operatorname{disc}}

\numberwithin{equation}{section}

\begin{document}

\title{The distribution of spacings between quadratic residues}
\author{P\"ar Kurlberg and Ze\'ev Rudnick} \address{Raymond and
  Beverly Sackler School of Mathematical Sciences, Tel Aviv
  University, Tel Aviv 69978, Israel} \date{Dec 14, 1998}
\thanks{Supported in part by a grant from the Israel Science
  Foundation. In addition, the first author was partially supported by
  the EC TMR network "Algebraic Lie Representations", EC-contract no
  ERB FMRX-CT97-0100}

\begin{abstract}

We study the distribution of spacings between squares modulo $q$, 
where $q$  is square-free and highly composite, in the limit as the 
number of prime  factors of $q$ goes to infinity. We show that all
correlation functions are Poissonian, which among other things, 
implies that the spacings between nearest neighbors, normalized to have unit
mean, have an exponential distribution.

\end{abstract}

\maketitle




\section{Introduction}

Our goal in this paper is to study the distribution of spacings (or
gaps) between squares in $\Z/q\Z$, as $q\to \infty$.
In the case that $q$ is {\em prime}, 
a theorem of Davenport 
\cite{Davenport-quadratic, Davenport-powers, Katz, Schmidt-book} 
shows that the  
probability of two consecutive quadratic residues modulo a prime $q$
being spaced $h$ units apart is $2^{-h}$, as $q\to\infty$. 
For our purposes, we may interpret this result as saying that when 
we normalize the spacings to have unit mean, 
then the distribution of spacing as $q\to \infty$ along primes is
given by 
$$
P(s) = \sum_{h=1}^\infty 2^{-h} \delta(s-\frac h2)
$$
that is, a sum of point masses at half-integers with exponentially
decreasing weights. 

In this paper we study the spacing distribution of squares 
modulo $q$ when $q$ is  square-free and {\em highly composite}, that is
the limiting  
distribution of spacings between the squares modulo $q$ as the number
of prime divisors, $\omega(q)$, tends to infinity.  
For odd square-free $q$ the number $N_q$ of squares modulo $q$ equals 
$$
N_q=\prod_{p\mid q} \frac {p+1}2
$$
This is because if $p$ is an odd prime,  
the number of squares modulo $p$ is $(p+1)/2$ and 
for $q$ square-free, $x$ is a square modulo $q$ if and only if
$x$ is a square modulo $p$ for all primes $p$ dividing $q$.  
Thus for odd $q$, the mean spacing $s_q=q/N$ equals 
$$
s_q=\frac{2^{\omega(q)}}{\prod_{p\mid q} (1+1/p)} =
\frac{2^{\omega(q)}}{\sigma_{-1}(q)}  
$$
For $q=2q'$ even and square-free, it is easily seen that $s_{q}=s_{q'}$. 
It follows that $s_q \to \infty$ as
$\omega(q)\to \infty$, unlike the case of prime $q$ where the mean
spacing is essentially constant. Thus, unlike in the prime case where 
the level spacing distribution was forced to be supported on a lattice, 
in the highly composite case there is an a-priori chance of getting a
continuous distribution. 

A relevant statistical model  for the distribution of
spacings is given by 
looking at {\em random} points in the unit interval $\R/\Z$. 
For independent, uniformly distributed 
numbers in $\R/\Z$,  the spacing statistics are said to be 
{\em  Poissonian}. The distribution $P(s)$ of spacings between
consecutive points will be that of a Poisson arrival process, i.e.
$P(s) = e^{-s}$ (see \cite{Feller}). 
Moreover, the joint distribution of $k$ consecutive
spacings is the product of $k$ independent exponential random
variables.

It is well known \cite{Mehta} that the spacing statistics  of the
superposition of several {\em independent} spectra converges to the
Poisson case - the spacings statistics of {\em uncorrelated}  levels. 
Thus the heuristic that ``primes are independent'' together with 
Davenport's result  indicates that the spacing statistics of the
squares modulo $q$ 
should  in the limit as $\omega(q)\to \infty$  be Poissonian, 
i.e., that in some sense
squares modulo $q$ behaves as random numbers. It is our purpose to
confirm this expectation. 


In order to study the level spacings, we proceed by studying the 
{\em   $r$-~level correlation functions}. These  measure clustering
properties of a sequence in $\R/\Z$ on a scale of the mean spacing.   
Their  definition and their application to computing various local
spacings statistics are recalled in appendix~\ref{app:spacings}. 
In our case, these turn out to be given by  by the following: 
For $r\geq2$ and a bounded convex set $\CC\subset \R^{r-1}$, let
$$
R_r(\CC,q) = \frac 1{N_q}\#\{x_i \mbox{ distinct squares mod } q : 
(x_1-x_2,\dots x_{r-1}-x_r) \in s\CC \}.
$$
This is immediately transformed into 
\begin{equation}
R_r(\CC,q)  = \frac 1N_q \sum_{h\in s\CC\cap \Z^{r-1}} N(h,q)
\end{equation}
where $N(h,q)$ is the number of solutions of the system of congruences 
$y_{i+1}-y_i=h_i \mmod q$ with $y_1, y_2, \ldots y_r$    
squares modulo $q$ and $h = (h_1, \ldots h_{r-1}) \in \Z^{r-1}$. 

To compute the correlations for {\em distinct} $x_i$ we consider only
sets $\CC$  which a-priori only contain vectors $(x_i-x_{i+1})$ with 
distinct coordinates. To do  this, we define ``roots'' 
$\sigma_{ij}$ on $\R^{r-1}$ for $i<j$ by $\sigma_{ij}(h) =
\sum_{k=i}^{j-1}h_k$. The hyper-planes $\{\sigma_{ij}=0\}\subset
\R^{r-1}$ are called ``walls'', and $(x_i-x_{i+1})$ does not lie in
any of the walls if and only if all coordinates $x_i$ are distinct.

Our main result 
shows that if $\CC$ does not intersect any wall then $R_r(\CC,q) \to
\vol(\CC)$ for any sequence of  square-free $q$ with
$\omega(q)\to\infty$:  
\begin{thm}\label{rlevel thm}
Let $q$ be  square-free,  $r\geq 2$ and
$\CC\subset \R^{r-1}$ a bounded convex set which does not intersect
any of the walls. Then 
the $r$-~level correlation function satisfies 
$$
R_r(\CC,q) = \vol(\CC) + O(s^{-1/2+\epsilon}) \qquad \mbox{as } s\to \infty
$$
for all $\epsilon>0$, where $s$ is the mean spacing. 
\end{thm}

This theorem  implies that all spacing statistics are Poissonian 
(see Appendix~\ref{app:spacings}). 
For instance, if we denote  by $s_1,\dots, s_{N-1}$ the normalized
differences between neighboring squares, then 
we have 
\begin{thm} \label{thm: spacings}
For $q$  square-free, the limiting level spacing distribution of
the squares modulo $q$ is given by $P(t) =\exp(-t)$ as
$\omega(q)\to\infty$. Moreover, under the same condition, 
for any $k\geq1$ the limiting joint distribution of
$(s_n,s_{n+1},\dots,s_{n+k})$ is  a product $\prod_{i=0}^k\exp(-t_i)$
of $k+1$ independent exponential variables. 
\end{thm}

There are only a few known cases  where the complete spacing distribution
can be proved to be Poissonian as in our case. A  notable example is 
Hooley's results \cite{Hooley1, Hooley2, Hooley3, Hooley-survey} that the 
spacings between elements co-prime to $q$  are 
Poissonian as the mean spacing $q/\phi(q)\to\infty$. 
A much more recent result is due to Cobeli and Zaharescu \cite{CZ}
who show that the spacings between primitive roots modulo a prime $p$
are Poissonian provided the mean spacing $p/\phi(p-1)\to \infty$. 

The results of this paper are related to work on the level spacing
distribution of the fractional parts $\{\alpha n^2\}$ ($\alpha$
irrational) by Rudnick, Sarnak and Zaharescu
\cite{Rudnick-Sarnak,RSZ}. 
In particular, in \cite{RSZ} an attempt to study that
problem is made by replacing   $\alpha$ with a rational approximation
$b/q$, and this leads to study the spacings of the sequence $b n^2
\mod q$, $1\leq n\leq N$ for $N$ a small power of $q$. 
The available {\em sites} are exactly the set of squares modulo $q$,
and hence our interest in the problem. 

In \cite{RSZ}, it is shown that in order for all the correlation
functions of the sequence $\{\alpha n^2\}$ to have Poisson behavior,
it is necessary to assume that the rational approximants $b/q$ have
denominator $q$ which is close to square-free. Hence our interest 
in the square-free case. 
For {\em arbitrary} $q$ it is still true that all correlations are
Poissonian, but there are significant technical complications to
overcome in proving this, see \cite{kurlberg}. 

We believe that the
methods developed in this paper should be useful in studying similar
problems, for instance the spacing distribution of cubes modulo $q$,
as the number of prime factors of $q$ that are congruent to $1$ modulo
$3$ tends to infinity.  (The condition modulo $3$ is necessary in
order for the mean spacing to go to infinity.)

\subsection*{Contents of the paper:} 
We begin with a section sketching the argument for 
Theorem~\ref{rlevel thm} in the case of the pair correlation
function. This section can be used as a guide to the rest of the
paper.


In section~\ref{sec:odd} we first reduce the problem to the case that
$q$ is {\em odd}.  
Then in section~\ref{sec:primecase} we analyze the behavior of
$N(h,p)$ where $p$ is prime. Squares that are distinct
modulo $q$ are not necessarily distinct modulo $p$; we denote by 
$r_{eff}(h)$  the number of squares that remains distinct after
reduction modulo $p$. Using an inclusion-exclusion argument we write
$r_{eff}(h)$ as a linear combination of characteristic functions of
certain hyper-planes over $\Z/p\Z$.
Next, in section~\ref{sec:expression for R} we use the multiplicative 
properties of the counting functions $N(h,q)$  
to derive an expression for $R_r(\CC,q)$ as a sum
over divisors $c$ of $q$ and lattices $L$ arising from intersections
of hyper-planes modulo $p$ for different $p$'s (proposition
\ref{prop:Expression for R_r}). 

In section~\ref{sec:rlevel} we show that the main term of the 
sum consists of those terms for which the product of $c$ and the
discriminant of $L$ are small with respect to $s$, and an error term
corresponding to terms where the product is large. 
In section~\ref{sec:mainterm} we evaluate
the main term 
and show that it gives us exactly $\vol(\CC)$ , thus giving us our main
result. 

In  appendix~\ref{app:spacings}  we explain how to use 
Theorem~\ref{rlevel thm} to derive results such as 
Theorem~\ref{thm: spacings}, 
that the level spacings are Poissonian as well. 
Appendix~\ref{sec:geom} explains some background on counting
lattice points in  convex sets. 
In appendix~\ref{app:divisors} we 
estimate  the number of divisors of $q$ that are smaller than 
a fixed power of the mean spacing $s$.

\newpage
\section{The pair correlation - a sketch}

In order to explain the proof of our main theorem \ref{rlevel thm}, 
we give an overview of the argument in the special 
case of the pair correlation function. 

Let $q$ be an odd, square-free number with $\omega(q)$ prime factors,
and $I$  an interval, not containing the origin. Define as in the 
introduction  the pair correlation function 
$$
R_2(I,q) = \frac 1N \sum_{h\in sI\cap \Z} N(h,q)
$$
where $N$ is the number of squares modulo $q$, 
$s=q/N=2^{\omega(q)}/\sigma_{-1}(q)$ is their mean spacing, 
$\sigma_{-1}(q) = \prod_{p\mid q}(1+\frac 1p)$,  
and $N(h,q)$
is the number of solutions in squares modulo $q$ of the equation 
$$
y_1-y_2=h \mod q
$$
We will sketch a proof that $R_2(I,q)\to |I|$ as
$\omega(q)\to\infty$ ($|I|$ being the length of the interval). 
In fact we have the more precise result
\begin{thm}
For $q$ odd, square-free we have for all $\epsilon>0$ 
$$
R_2(I,q) = |I|+O(s^{-1+\epsilon})
$$
\end{thm}
Here are the main steps in the argument: 

\noindent{\bf Step 1}: 

By the Chinese Remainder Theorem, $N(h,q) = \prod_{p\mid q}N(h,p)$ is
a product over primes dividing $q$. 
By elementary considerations, one sees that 
\begin{equation}\label{pc:N(h,p)}
N(h,p) = \frac{p+a(h,p)}4\Delta(h,p)
\end{equation}
with $a(h,p)=O(1)$ and 
$$
\Delta(h,p) = 1+\delta(h,p),
\qquad \delta(h,p) = \begin{cases} 0& p\nmid h\\1 &p\mid h \end{cases}
$$
 From this we see that 
\begin{equation}\label{pc:N(h,q)}
N(h,q) = \frac{q\Delta(h,q)}{4^{\omega(q)}} \sum_{c\mid q}\frac{a(h,c)}c
\end{equation}
with $a(h,c):=\prod_{p\mid c}a(h,p) \ll c^{\epsilon}$ and 
$\Delta(h,q) = \prod_{p\mid q}\Delta(h,p)$.

\noindent{\bf Step 2:}

We decompose $\Delta(h,q) = \Delta(h,c)\Delta(h,\frac qc)$ and rewrite 
$\Delta(h,\frac qc)$ as 
$$
\Delta(h,\frac qc)=\prod_{p\mid \frac qc} (1+\delta(h,p)) =
\sum_{g\mid \frac qc} \delta(h,g)
$$
with 
$$
\delta(h,g)=\begin{cases}0&g\nmid h\\  1&g\mid h \end{cases}
$$
Substituting this into  the expression \eqref{pc:N(h,q)} for $N(h,q)$
and inserting the result into the
formula for $R_2(I,q)$, we get 
\begin{equation}\label{pc:R_2}
R_2(I,q) = \frac 1{\sigma_{-1}(q)2^{\omega(q)}} \sum_{c\mid q} \frac 1c 
\sum_{g\mid \frac qc} \sum_{h\in sI \cap g\Z}a(h,c)\Delta(h,c)
\end{equation}

\noindent{\bf Step 3:} 

We partition the sum into two parts, one over the pairs $g,c$ with
$gc<s$ and the leftover part over pairs with $gc\geq s$. 
We will show this leftover part is negligible (in fact
$O(s^{-1+\epsilon})$):  
We first use $a(h,c)\Delta(h,c)\ll c^\epsilon$ and the fact that in order
for the inner sum over $h$ to be nonempty, we need $g\ll s$ (recall
that $I$ does not contain the origin!) to get that the sum over pairs
with $cg>s$ is bounded by 
\begin{equation*}
\begin{split}
s^{-1+\epsilon}\sum_{c\mid q} c^{-1+\epsilon} 
\sum_{\substack{g\mid \frac qc\\ g\ll s\\cg>s}} \#(sI\cap g\Z)  
&\ll s^{-1+\epsilon}\sum_{c\mid q} c^{-1+\epsilon} 
\sum_{\substack{g\mid \frac qc\\ g\ll s\\cg>s}} \frac sg \\
&\ll s^{\epsilon}\sum_{\substack{d\mid q\\d>s}} d^{-1+\epsilon} 
\sum_{\substack{g\mid d \\ g\ll s}} 1 
\end{split}
\end{equation*}
Now we use Lemma \ref{small divisors} which shows that the number of
divisors $g<s$ of $q$ is a most $O(s^\epsilon)$ and
Lemma~\ref{bounding divisor sums} to bound the above by   
$$
s^{\epsilon}\sum_{\substack{d\mid q\\d>s}} d^{-1+\epsilon} \ll
s^{-1+\epsilon} 
$$
as promised.
 
\noindent{\bf Step 4:} 

For each pair of $c,g$ with $cg<s$, we first
treat the inner sum over $h\in sI\cap g\Z$. We break it up into sums
over $\frac{s|I|}{gc}+O(1)$ subintervals $[y,y+cg)\cap g\Z$ plus a
leftover term of size at most 
$c^{1+\epsilon}$. For each subinterval, we use periodicity of
$a(h,c)\Delta(h,c)$ under $h\mapsto h+c$ to find 
$$
\sum_{h\in [y,y+cg)\cap g\Z} a(h,c)\Delta(h,c) = 
\sum_{h_1=1}^c a(gh_1,c)\Delta(gh_1,c)
$$

Because $q$ is square-free, and $g$ divides $q/c$, we have that $g,c$
are coprime. Therefore we can change variables $h=gh_1$ to get that
this last sum equals 
$$
\sum_{h\mod c} a(h,c)\Delta(h,c) = \prod_{p\mid c} 
\sum_{h\mod p}a(h,p)\Delta(h,p)
$$

We evaluate the sum $\sum_{h\mod p}a(h,p)\Delta(h,p)$ by noting that 
summing \eqref{pc:N(h,p)} over $h \mod p$, the sum of the LHS is
simply the number of all pairs of squares modulo $p$, namely
$(p+1)^2/4$. This gives 
$$
\sum_{h\mod p}a(h,p)\Delta(h,p) = p+1
$$
Thus the inner sum over $h\in sI\cap g\Z$ equals 
\begin{equation*}
\begin{split}
\sum_{h\in sI \cap g\Z}a(h,c)\Delta(h,c) &=
\left(\frac{s|I|}{gc}+O(1) \right)\prod_{p\mid c} (p+1)+O(c^{1+\epsilon}) \\
&= \frac{s|I|}{g}\sigma_{-1}(c)+O(c^{1+\epsilon})
\end{split}
\end{equation*}

\noindent{\bf Step 5:}

Inserting this into the expression \eqref{pc:R_2} for $R_2(I,q)$ gives 
$$
R_2(I,q) = \frac 1{2^{\omega(q)} \sigma_{-1}(q)} \sum_{c\mid q} \frac 1c
\sum_{g\mid \frac qc: gc<s} \frac{s|I|}{g}\sigma_{-1}(c)
+O(s^{-1+\epsilon}) 
$$

Now we extend the sum to all pairs $g,c$,  to
find that up to an error of $O(s^{-1+\epsilon})$ we have 
\begin{equation*}
\begin{split}
R_2(I,q) &\sim |I| \frac 1{\sigma_{-1}(q)^2} \sum_{c\mid q} \frac
{\sigma_{-1}(c)}c \sum_{g\mid \frac qc} \frac 1g \\ 
&=  |I| \frac 1{\sigma_{-1}(q)^2} \sum_{c\mid q}
\frac{\sigma_{-1}(c)}c \sigma_{-1}(\frac qc) \\
&= |I| \frac 1{\sigma_{-1}(q)} \sum_{c\mid q}\frac 1c = |I|
\end{split}
\end{equation*}
which is what we need to prove our theorem \qed. 

In the following sections, we will repeat these steps with full
details for the higher
correlation functions, where several technical complications arise. 

\newpage
\section{Reduction to odd $q$} \label{sec:odd}

We first show that in Theorem~\ref{rlevel thm} it suffices to consider
only the case of $q$ {\em odd}:    
Suppose that $q=2q'$ with $q'$ odd and square-free. 
We recall that 
\begin{equation}\label{trans Rr}
R_r(\CC,q)  = \frac 1N_{q} \sum_{h\in s\CC\cap \Z^{r-1}} N(h,q)
\end{equation}
where $N(h,q)$ is the number of solutions of the system
$y_{i+1}-y_i=h_i$ where $y_1, y_2, \ldots y_r$    
are squares modulo $q$ and $h = (h_1, \ldots h_{r-1}) \in
({\Z}/q{\Z})^{r-1}$.

By the Chinese Remainder Theorem, the number $N_{q}$ of squares modulo
$q$ is the product  
$$
N_{q}=N_2 N_{q'} = 2N_{q'}
$$
Therefore the mean spacing $s_{q}:=q/N_q$ is given by 
\begin{equation}\label{newspacing}
s_q = \frac{2q'}{2N_{q'}} =  \frac{q'}{N_{q'}} =s_{q'}
\end{equation}
Moreover, again by the Chinese Remainder Theorem, 
$$
N(h,q) = N(h,2)N(h,q') 
$$
and since all residues modulo $2$ are squares, we have $N(h,2)=2$. 
Thus we find 
\begin{equation}\label{newratio}
\frac{N(h,q)}{N_{q}} = \frac{2N(h,q')}{2N_{q'}} = \frac{N(h,q')}{N_{q'}}
\end{equation}
Inserting \eqref{newspacing}, \eqref{newratio}  into \eqref{trans Rr},
we find that  
$$
R_r(\CC,q)   = R_r(\CC,q')  
$$
This shows that it suffices to prove Theorem~\ref{rlevel thm}  for $q$
odd, which we assume is the case in the sequel.   

\newpage
\newpage
\section{The prime case}\label{sec:primecase}

Let $p>2$ be a prime. 
For $h=(h_1,\dots h_{r-1}) \in (\Z/p\Z)^{r-1}$, we define $N_r(h,p)$ to be 
the number of solutions in {\em squares} $y_i$ mod $p$ (including $y_i=0$) 
of the system 
\begin{equation}\label{system}
y_i- y_{i+1} = h_i \mod p,\qquad 1\leq i \leq r-1
\end{equation}
This number depends crucially on the number of {\em distinct} $y_j$. 
For each $h= (h_1,\dots ,h_{r-1})$, we define 
$\reff(h)$ to be the number of distinct $y_j$ (not necessarily squares) 
satisfying the system (\ref{system}). 
Since the solutions of the {\em homogeneous} system $y_i- y_{i+1} = 0 \mod p$ 
are spanned by $(1,\dots,1)$, $\reff(h)$ is well-defined 
(independent of the particular solution $y$ of (\ref{system})). 

We define {\em roots} $\sigma_{ij}(h)$, $1\leq i<j\leq r$ by 
\begin{equation}
\sigma_{ij}(h) = \sum_{k=i}^{j-1} h_k
\end{equation}
so that $\sigma_{i,i+1}(h) = h_i$, 
$\sigma_{ij}= \sum_{k=i}^{j-1}\sigma_{k,k+1}$. 
The solutions of (\ref{system}) are all distinct of and only if 
$\sigma_{ij}(h) \neq 0$, for all $i<j$, since 
$$
y_i-y_j = 
\sum_{k=i}^{j-1} y_k-y_{k+1} = \sum_{k=i}^{j-1} h_k =\sigma_{ij}(h)
$$

\begin{prop}\label{propN(h,p)}
Let $\reff(h)$ be the number of {\em distinct} $y_i$ in a solution of 
(\ref{system}). Then 
\begin{equation}\label{formula for N(h,p)}
N_r(h,p) = \frac{p+a(h,p)}{2^{\reff}}
\end{equation}
with $a(h,p)\ll_r p^{1/2}$. 
\end{prop}
\begin{proof}
The case $\reff(h)=1$ is precisely when $h=0$ and all $y_i$ are equal:
$y_1=y_2=\dots = y_r$. In this case the number of solutions is the
number of squares modulo $p$, namely $(p+1)/2$, which is of the desired
form. We thus assume from now that $\reff(h)>1$. 

We first reduce the system (\ref{system}) to a system of $\reff-1$ 
equations in $\reff$ variables: If $\reff(h)$ is the number of distinct 
$y_i$ in a solution of (\ref{system}) (independent of $y$!), then 
we can eliminate some of the equations. Renumber the variables so that 
$y_1,\dots y_{\reff}$ 
are the {\em distinct} coordinates of a solution, and for all 
$j\geq 1$, $y_{\reff+j}$ 
equals one of these, then the system  (\ref{system}) is equivalent to the 
reduced system 
\begin{equation}\label{systemred}
y_i-y_{i+1} = h'_i \mod p, \qquad 1\leq i\leq \reff-1
\end{equation}
(where the $h'_i$ are renumbered $h_j$ to give that the first $\reff$ 
coordinates  are distinct).
So we need to find the number of solution of the reduced system 
(\ref{systemred}).

We first eliminate those solutions where at least one of the $y_j$ is zero. 
In this case, since the system (\ref{systemred}) (considered as 
a {\em linear} system) has rank $\reff-1$ in  $\reff$ variables, specifying 
any one of the variables determines all the others, hence the number 
of solutions with some coordinate zero is at most $\reff$. Thus we need only 
count solutions where all coordinates $y_i$ are nonzero. 

To every such solution in squares $y_i\neq 0 \mod p$, write $y_i=x_i^2 \mod p$ 
with $x_i\neq 0 \mod p$. There are precisely two such solutions, 
namely $\pm x_i \mod p$. Thus the number of possible $x_i$ 
corresponding to a given solution $y$ of (\ref{systemred}) is precisely 
$2^{\reff}$, and the number of nonzero solutions of the reduced system 
(\ref{systemred}) with $y_i$ squares  modulo $p$  is exactly 
$1/2^{\reff}$ times the number of solutions of the system 
\begin{equation}
x_i^2-x_{i+1}^2 = h'_i \mod p, \qquad 
1\leq i \leq \reff-1
\end{equation}
with $x_i \neq 0 \mod p$. By adding back at most $r$ solutions we can 
remove the condition $x_i \neq 0$, and then we find that 
\begin{equation}\label{N in terms of n}
N_r(h,p) =  \frac 1{2^{\reff}} n(h',p) + O_r(1)
\end{equation}
where $n(h',p)$ is the number of solutions of 
$$
x_i^2-x_{i+1}^2 = h'_i \mod p\qquad 
1\leq i \leq \reff(h)-1
$$
This is just the number of solutions $(t,x_1,\dots ,x_{\reff})$ 
of the system
\begin{equation}
x_1^2 = t-b_1,\quad x_2^2=t-b_2,\dots, x_{ \reff}^2 = t-b_{\reff} 
\end{equation}
with $b_1=0$, $b_2=h'_1$, $b_3=h'_1+h'_2$, \dots, 
$ b_{\reff(h)}  =h'_1+h'_2+\dots + h'_{\reff-1}$ and in general 
$b_k=\sigma_{1k}(h')$. 
Note that the $b_i$ are {\em distinct} 
 - this is equivalent to the requirement that the  solutions of the
reduced system~\eqref{systemred} be distinct. 
One can  now use the  ``Riemann Hypothesis for curves'' \cite{Weil} 
(see Schmidt's book \cite{Schmidt}, Chapter II, Theorem 5A and
Corollary 5B for the case $b_1=-1,b_2=-2,\dots b_r=-r$), 
to find 
\begin{equation}
|n(h',p) - p| \ll \reff 2^{\reff}\sqrt{p}
\end{equation}

In addition, $|N(h,p) - n(h',p)/2^{\reff}|\leq r$ and so 
$$
N(h,p) = \frac{p+a(h,p)}{2^{\reff}}
$$
with 
$$
a(h,p) \ll 2^{\reff}(\reff \sqrt{p}+r) \ll_r \sqrt{p}
$$
This proves Proposition \ref{propN(h,p)}. 
\end{proof}


\subsection{A formula for $\reff(h)$} 
Our next order of business is to give a formula  for $\reff(h)$. 
We begin with some combinatorial background: 
A {\em set partition} of the set $\{1,2,\dots,r\}$ is a collection of 
disjoint subsets $\sF = [F_1,\dots, F_t]$, $F_i\subseteq \{1,2,\dots,r\}$, 
whose union is all of $\{1,2,\dots,r\}$. We set $|\sF| = t$, the number of 
subsets in $\sF$. 

To each set partition $\sF$, we associate a subset $V_{\sF}$ of affine 
$r$-space $V=\Af^r$ by setting 
\begin{equation}
V_{\sF} = \{ s\in \Af^r: s_i = s_j \mbox{ if } i,j \mbox{ are in some } F_k\}
\end{equation}
Correspondingly, in $H=\Af^{r-1}$ 
we have a  subspace 
\begin{equation}
H_{\sF} = \{ h\in \Af^{r-1}:\sigma_{ij}(h) = 0\mbox{ if } 
i,j \mbox{ are in some } F_k\}
\end{equation}
Under the map $\pi:V\to H$ taking $s=(s_i)\mapsto (s_i-s_{i+1})$, 
we have $V_{\sF} = \pi^{-1}H_{\sF}$. 

There is a partial ordering on the collection of all set-partitions of 
$\{1,\dots,r\}$ with $\sF \preceq \sG$ if and only if  every $F_i$ is contained in 
some $G_j$. 

For example, $\sO = [\{1,2,\dots,r\}]$ is the maximal element of this 
partial ordering, with $|\sO|=1$ and  $H_{\sO} =  (0)$. 
The minimal element is 
$\sR = [\{1\},\{2\},\dots,\{r\}]$ with $|\sR| = r$ and $H_{\sR} = \Af^{r-1}$.

The partial ordering on set-partitions is inclusion-reversing on subspaces: 
$\sF\preceq \sG\Leftrightarrow V_{\sF} \supseteq V_{\sG} \Leftrightarrow 
 H_{\sF} \supseteq H_{\sG}$. 

The {\em regular} part of $V_{\sF}$ is 
$$
V_{\sF}^\times = \{s\in V_{\sF}: s_i\neq s_j \mbox{ if } i,j 
\mbox{ are not in some } F_k\}
$$
and likewise we define $H_{\sF}^\times$.  
Then $H_{\sO}^\times = H_{\sO} = (0)$, and every $h$ belongs to a unique 
$H_{\sF}^\times$ for some $\sF$. We thus have 
$$ 
H= \coprod_{\sF} H_{\sF}^\times
$$
and likewise
$$
H_{\sF}= \coprod_{\sF\preceq \sG} H_{\sG}^\times
$$

We can now give a formula for $\reff(h)$: 
\begin{equation} 
\reff(h)  = \dim V_{\sF} = \dim H_{\sF} + 1 = |\sF|
\end{equation}
where $\sF$ is the unique set-partition such that $h\in H_{\sF}^\times$. 

We can write this as follows: Define 
\begin{equation}
\delta_{\sF}(h) = \begin{cases} 
                      1& h\in H_{\sF} \\
                      0&\mbox{otherwise}
                   \end{cases}
,\qquad 
\delta_{\sF}^\times(h)= \begin{cases} 
                      1& h\in H_{\sF}^\times \\
                      0&\mbox{otherwise}
                   \end{cases}
\end{equation}
Then 
\begin{equation}
\reff(h) = \sum_{\sF}\dim (V_{\sF}) \delta_{\sF}^\times(h)
\end{equation}
Similarly 
\begin{equation}\label{def Delta}
\Delta(h,p): = 2^{r-\reff(h)} = \sum_{\sF}
2^{\codim (V_{\sF})} \delta_{\sF}^\times(h)
\end{equation}
It will be convenient to express this in terms of the characteristic 
function $\delta_{\sF}$ of the {\em subspaces} $H_{\sF}$. For this we use 
M\"obius inversion. Since the collection of all set-partitions of 
$\{1,\dots,r\}$ is a partially-ordered set, it has a M\"obius function 
$\mu(\sF,\sG)$ which is the unique function so that for any functions 
$\psi$, $\phi$ on set-partitions satisfying 
\begin{equation} 
\phi(\sF) = \sum_{\sF\preceq \sG} \psi(\sG)
\end{equation}
we have 
\begin{equation}
%
%
\psi(\sF) = \sum_{\sF\preceq \sG} \mu(\sF,\sG)  \phi(\sG)
\end{equation}
An explicit form of $\mu(\sF,\sG)$ can be found in \cite{vLW}, \S 25. 
We will not have any use for it.

In our case, clearly we have $H_{\sF} = 
\coprod_{\sF\preceq \sG} H_{\sG}^\times$ so that 
\begin{equation}
\delta_{\sF} = \sum_{\sF\preceq \sG} \delta_{\sG}^\times
\end{equation}
Thus we have
\begin{equation}\label{deltasFtimes}
\delta_{\sF}^\times = \sum_{\sF\preceq \sG} \mu(\sF,\sG)\delta_{\sG}
\end{equation}
This gives us the formula for $\Delta(h,p) = 2^{r-\reff(h)}$: From 
\eqref{def Delta} and \eqref{deltasFtimes} we find
\begin{equation}\label{Delta(h,p)}
\Delta(h,p) = \sum_{\sG} \lambda(\sG) \delta_{\sG}(h)
\end{equation}
with 
\begin{equation}\label{lambda(sG)}
\lambda(\sG)=\sum_{\sF\preceq \sG} \mu(\sF,\sG)2^{\codim V_{\sF}}
\end{equation}

For use in Section~\ref{sec:mainterm}, we need to know the sum of the
product of $\Delta(h,p)$ with the error term $a(h,p)$ in 
\eqref{formula for N(h,p)} over all vectors $h$:  
\begin{lem}\label{sum a(h,p)Delta(h,p)}
$$
\sum_{h\mod p} a(h,p)\Delta(h,p) = 
(p+1)^r - p^r\sum_{\sG}\lambda(\sG) p^{-\codim H_{\sG}}
$$
\end{lem}
\begin{proof}
We have by definition  
$$
N(h,p) = \frac{p+a(h,p)}{2^r} \Delta(h,p)
$$
so that 
$$
a(h,p)\Delta(h,p) = 2^r N(h,p) - p \Delta(h,p)
$$
Now sum over all $h \mod p$: The sum of $N(h,p)$ is just the total
number of $r$-tuples of squares modulo $p$, namely $(\frac{p+1}2)^r$. 
To sum $\Delta(h,p)$ over $h$, we use \eqref{Delta(h,p)}: 
Since the sum over all $h$ of $\delta_{\sG}(h)$ is just the number of
vectors in the subspace $H_{\sG}$, namely 
$p^{\dim  H_{\sG}} =p^{r-1-\codim H_{\sG}}$, we find 
\begin{equation*}
\begin{split}
\sum_{h\mod p} a(h,p)\Delta(h,p) &= (p+1)^r - 
p\sum_{h\mod p} \Delta(h,p)  \\
&= (p+1)^r - p^r\sum_{\sG} \lambda(\sG)
p^{-\codim H_{\sG}} 
\end{split}
\end{equation*}
as required. 
\end{proof}

\newpage

\section{A formula for $R_r(\CC,q)$} \label{sec:expression for R}

In order to prove Theorem~\ref{rlevel thm}, we give an expression
(\ref{expr for R_r}) for  
the $r$-level correlation $R_r(\CC,q)$ which involves summing over 
the intersection of the dilated set $s\CC$ with various lattices. 

Recall that for each set-partition $\sG$ of $\{1,\dots,r\}$ we
associated a subspace $H_{\sG}\subseteq (\Z/p\Z)^{r-1}$. 
Now given a divisor $d\mod q$, let $\sGG=\otimes_{p\mid d}\sGp$,  
be a tuple of such set-partitions, one for each prime $p$ dividing
$d$ (recall that $q$, hence $d$, is square-free). 
Let $L(\sGG)\subset \Z^{r-1}$ be the pre-image of $\prod_{p\mid d}
H_{\sGp}$ under the reduction map $\Z^{r-1} \to \prod_{p\mid d}
(\Z/p\Z)^{r-1} \simeq (\Z/d\Z)^{r-1}$. 
$L(\sGG)$ is a lattice,  whose discriminant (that  is, the index in
$\Z^{r-1}$) is     
$$
\disc(\sGG)   = \prod_{p\mid d} p^{\codim(H_{\sGp})}
$$
The {\em support}  $\supp(\sGG)$ of $L(\sGG)$ is the product of all
primes $p$ for which $H_{\sGp}\neq (\Z/p\Z)^{r-1}$:  
$$
\supp(\sGG) = \prod_{p: \sGp \neq [\{1\},\dots,\{r-1\}]}p
$$
Since $\codim (H_{\sGp}) \leq r-1$, we get 
$$
\supp(\sGG) \mid \disc(\sGG) \mid \supp(\sGG)^{r-1}
$$

We  set
$$
\lambda(\sGG) = \prod_{p\mid d} \lambda(\sGp)
$$
where $\lambda(\sG)$ is given by \eqref{lambda(sG)}. 
We also set for a divisor $c\mid q$
$$
a(h,c):=\prod_{p\mid c}a(h,p), \qquad \Delta(h,c) := \prod_{p\mid c}\Delta(h,p)
$$
Note that by Proposition \ref{propN(h,p)}
\begin{equation}\label{bound for a(h,c)}
a(h,c)\ll c^{1/2+\epsilon}, \qquad \Delta(h,c)\ll c^\epsilon
\end{equation}
for all $\epsilon>0$. 

Our formula for $R_r(\CC,q)$ is 
\begin{prop}\label{prop:Expression for R_r}
The $r$-level correlation function is given by
\begin{equation}\label{expr for R_r} 
R_r(\CC,q)=  \frac{s}{2^{r\omega(q)}} \sum_{c\mid q} \frac 1c
\sum_{\supp(\sGG)\mid \frac qc} \lambda(\sGG) 
\sum_{h\in s\CC\cap L(\sGG)}a(h,c) \Delta(h,c)
\end{equation}
\end{prop}
\begin{proof}
We have that 
\begin{equation*}
R_r(\CC,q)  = \frac 1N \sum_{h\in s\CC\cap \Z^{r-1}} N(h,q)
\end{equation*}
By the Chinese Remainder Theorem, 
$$
N(h,q) = \prod_{p\mid q} N(h,p)
$$
We rewrite formula (\ref{formula for N(h,p)}) in the form 
$$
N(h,p) = \frac{p+a(h,p)}{2^r} \Delta(h,p)
$$
where 
$$
\Delta(h,p) = 2^{r-\reff(h)}
$$
Thus we find 
\begin{equation}\label{formula for N(h,q)}
N(h,q) = \frac{q\Delta(h,q)}{2^{r\omega(q)}}
\sum_{c\mid q} \frac{a(h,c)}c = 
\frac{q}{2^{r\omega(q)}}
\sum_{c\mid q}\Delta(h,\frac qc) \frac{a(h,c)\Delta(h,c)}c
\end{equation}
Inserting \eqref{formula for N(h,q)} we get a formula for $R_r(\CC,q)$:
Recalling that $N=q/s$,
\begin{equation}
R_r(\CC,q) = \frac s{2^{r\omega(q)}} \sum_{c|q}
\frac{1}c \sum_{h\in s\CC} \Delta(h,\frac qc)a(h,c) \Delta(h,c)
\end{equation}

Next we use the expression (\ref{Delta(h,p)}) for $\Delta(h,p)$ to
write $\Delta(h,\frac qc) = \prod_{p|q/c}\Delta(h,p)$ in the form 
\begin{equation} 
\Delta(h,\frac qc) = \prod_{p|\frac qc} \sum_{\sGp} \lambda(\sGp)
\delta(h,\sGp)  = \sum_{\sGG=\otimes_{p|\frac qc} \sGp}
\lambda(\sGG)\delta(h,\sGG) 
\end{equation} 
where  the sum is over all tuples of set-partitions 
$\sGG=\otimes_{p\mid \frac qc}\sGp$,  one for each prime dividing
$\frac qc$, 
and  we put for each such tuple $\sGG$
$$
\lambda(\sGG) := \prod_{p|\frac qc} \lambda(\sGp)
$$
and   
$$
\delta(h,\sGG):= \prod_{p|\frac qc} \delta(h,\sGp) = 
\begin{cases} 1& h\in H_{\sGp} \mod p \mbox{ for all }p|\frac qc\\
              0&\mbox{otherwise} 
\end{cases} 
$$
This is the characteristic function of the lattice $L(\sGG)$ 
whose support $\supp(\sGG)$ divides $q/c$. 
Thus we get the desired expression for $R_r(\CC,q)$ 
\begin{equation*}
R_r(\CC,q)=  \frac{s}{2^{r\omega(q)}} \sum_{c\mid q} \frac 1c
\sum_{\supp(\sGG)\mid \frac qc} \lambda(\sGG) 
\sum_{h\in s\CC\cap L(\sGG)}a(h,c) \Delta(h,c)
\end{equation*}
\end{proof}


\newpage

\section{Evaluating the $r$-level correlations}\label{sec:rlevel}

In order to estimate the correlations using 
Proposition~\ref{prop:Expression for R_r},  we 
partition the sum \eqref{expr for R_r} into two parts: the first
consisting of pairs $c$ 
and $\sGG$ such that  $c\disc(\sGG) <s$, and the second of the pairs for which
$c\disc(\sGG)>s$. We will show that the first part gives the main term
and the second is negligible. 

\subsection{The case $c\disc(\sGG)>s$}\label{sec cgamma>s} 
We use $a(h,c)\ll c^{1/2+\epsilon}$ (\ref{bound for a(h,c)}), and 
$\Delta(h,c)\ll c^\epsilon$ to see that this term is bounded by 
\begin{equation}\label{terms with c disc(sGG)>s}
\frac{s}{2^{r\omega(q)}} \sum_{c\mid q} \frac 1c
\sum_{\substack{\supp(\sGG)\mid\frac qc\\c\disc(\sGG)>s}}
|\lambda(\sGG)| \#\{  s\CC\cap L(\sGG)\}c^{1/2+\epsilon}
\end{equation} 
By the Lipschitz principle (Lemma \ref{Lipschitz lemma}), 
$$
\#\{  s\CC\cap L(\sGG)\} \ll \frac{\vol(s\CC)}{\disc(\sGG)} + s^{r-2}
$$
and since $\vol(s\CC) = s^{r-1} \vol(\CC)$,  we find that 
\begin{equation}
 \#\{  s\CC\cap L(\sGG)\} \ll \frac{s^{r-1}}{\disc(\sGG)} + s^{r-2}
\end{equation}
Moreover, in order that $s\CC\cap L(\sGG) \neq \emptyset$, we will see
that we need $\supp(\sGG)
\ll s^{r(r-1)/2}$, since $\CC$ does not intersect the walls. 
This is a consequence of the following observation: 
Let $\CC\subset \R^{r-1}$ be a bounded convex set. 
Define 
$$
%
%
\diam_1(\CC) = \max\{\sum_{k=1}^{r-1} |x_k| : x\in \CC \}
$$
Note that $\diam_1$ scales linearly: $\diam_1(s\CC) = s\diam_1(\CC)$
for all $s>0$. 
\begin{lem} \label{avoid walls}
If $\supp(\sGG)> \diam_1(s\CC)^{r(r-1)/2}$ then $s\CC\cap
L(\sGG)$ is contained  in the walls $\{h\in \R^{r-1}: \sigma_{ij}(h) =
0 \mbox{ for some } i<j \}$. 
\end{lem} 
\begin{proof} 
Let $d_{ij}(\sGG)$ be the product of the primes $p$ such that
$\sigma_{ij}$ vanishes on $H_{\sGp}$, i.e. so that 
$$
\sigma_{ij}(x) = 0 \mod p \qquad \mbox{for all } x\in L(\sGG)
$$
Then $d_{ij}(\sGG)|\supp(\sGG)$ and moreover we claim that: 
$$
\disc(\sGG) | \prod_{i<j} d_{ij}(\sGG)
$$
It is enough to check this one prime at a time and is equivalent to
saying that 
$$
\codim(H_{\sGp}) \leq \#\{ i<j: \sigma_{ij} = 0 \mbox{ on } H_{\sGp}
\}
$$
which follows since $H_{\sGp}$ is given by vanishing of some of
the  $\sigma_{ij}$. 

Now note that if $\supp(\sGG)> d^{r(r-1)/2}$ then for some $i<j$, 
$d_{ij}(\sGG)>d$ because $\supp(\sGG)\leq \disc(\sGG) \leq\prod_{i<j}
d_{ij}(\sGG)$ and the last product consists of $r(r-1)/2$ factors.
If we  take $d=\diam_1(s\CC) = s\diam_1(\CC)$, then one has
$d_{ij}(\sGG)>\diam_1(s\CC)$ for some $i<j$. However $\sigma_{ij}(h) =
0 \mod d_{ij}(\sGG)$ and so 
$\sigma_{ij}(h) = m d_{ij}(\sGG)$ for some integer $m$. If $m=0$ then 
$h$ lies in a wall. If $m\neq 0$ then being an integer, $|m| \geq 1$
and so 
$$
|\sigma_{ij}(h)| \geq d_{ij}(\sGG)> \diam_1(s\CC)
$$
Since 
$$
\sigma_{ij}(h) = |\sum_{k=i}^{j-1} h_k| \leq \sum_{k=i}^{j-1} |h_k| 
\leq \sum_{k=i}^{r-1} |h_k| 
$$
we find that 
$$
\sum_{k=i}^{r-1} |h_k| > \diam_1(s\CC)
$$
Thus $h\not\in s\CC$ by definition of $\diam_1(s\CC)$.
\end{proof}
By Lemma \ref{avoid walls}, together with $|\lambda(\sGG)| \ll
\supp(\sGG)^\epsilon$, \eqref{terms with c disc(sGG)>s} is bounded by
\begin{equation}\label{truncated c sup>s}
\frac{s}{2^{r\omega(q)}} \sum_{c\mid q} c^{-1/2+\epsilon}
\sum_{\substack{\supp(\sGG)\mid \frac qc \\ c\disc(\sGG)>s \\ \supp(\sGG) \ll
s^{r(r-1)/2} }}
\supp(\sGG)^\epsilon
 \left(\frac{s^{r-1}}{\disc(\sGG)} + s^{r-2}\right)
\end{equation}
We split the sum into two parts and use $s< 2^{\omega(q)}$ to bound
\eqref{truncated c sup>s} by the sum of 
%
\begin{equation}
\label{truncated sum}
\frac 1s \sum_{c\mid q} c^{-1/2+\epsilon}
\sum_{\substack{\supp(\sGG)\mid \frac qc\\ c\disc(\sGG)>s  }}
\supp(\sGG)^\epsilon 
\frac{s}{\disc(\sGG)}
\end{equation}
and
\begin{equation}\label{extra term}
\frac 1s \sum_{c\mid q} c^{-1/2+\epsilon}
\sum_{\substack{\supp(\sGG)\mid \frac qc\\ c\disc(\sGG)>s \\ \supp(\sGG) \ll
s^{r(r-1)/2} }}
\supp(\sGG)^\epsilon 
\end{equation}

We begin by noting that the number of $\sGG$ with $\supp(\sGG) = g$ is
$O(g^\epsilon)$, i.e., 
\begin{equation}
\label{e:lattice-supp-sum}
\sum_{\supp(\sGG) = g} 1 \ll g^\epsilon,
\end{equation}

Since we sum over $\supp(\sGG) \ll s^{r(r-1)/2}$ in \eqref{extra term},
we have $\supp(\sGG)^\epsilon \ll s^{\epsilon'}$, and thus \eqref{extra
  term} is bounded by
$$
\frac 1s  \sum_{c\mid q} c^{-1/2+\epsilon} 
\sum_{\substack{ g\mid \frac qc \\g \ll s^{r(r-1)/2}}} g^{\epsilon}
\ll
s^{-1+\epsilon}  \sum_{c\mid q} c^{-1/2+\epsilon} 
\sum_{\substack{ g\mid \frac qc \\g \ll s^{r(r-1)/2}}} 1.
$$
By Lemma \ref{small divisors}, the number of divisors of $q/c$ which
are less than 
$s^{r(r-1)/2}$ is at most $s^\epsilon$, so this term is bounded by
$$
s^{-1+\epsilon}  \sum_{c\mid q} c^{-1/2+\epsilon}. 
$$
Since 
$$
\sum_{c\mid q} c^{-1/2+\epsilon}  = \prod_{p|q} (1+\frac
1{p^{1/2-\epsilon}}) \ll \prod_{p\mid q}(1+1)^{\epsilon'} \ll
s^{\epsilon''}
$$
the contribution of \eqref{extra term} is at most $O(s^{-1+\epsilon})$. 

It now remains to bound \eqref{truncated sum}. We first consider the
terms for which $c \supp(\sGG) > s$. Now, $\disc(\sGG) \geq
\supp(\sGG)$, so if $c\supp(\sGG) > s$ then certainly
$c\disc(\sGG) > s$, and sum of the corresponding terms in
\eqref{truncated sum} is bounded by
$$
\frac 1s \sum_{c\mid q} c^{-1/2+\epsilon}
\sum_{\substack{\supp(\sGG)\mid \frac qc\\ c\supp(\sGG)>s \\ }}
\supp(\sGG)^\epsilon \frac{s}{\supp(\sGG)} 
$$
$$
= \sum_{c\mid q} c^{-1/2+\epsilon}
\sum_{\substack{g \mid \frac qc\\cg>s }}
\frac{1}{g^{1-\epsilon}}
\sum_{\supp(\sGG) = g } 1
\ll
\sum_{c\mid q} c^{-1/2+\epsilon}
\sum_{\substack{ g\mid \frac qc \\ cg>s }} 
\frac 1{g^{1-\epsilon}},
$$
by \eqref{e:lattice-supp-sum}. Changing variable to $d=cg$, which is a
divisor of $q$ satisfying $d>s$, this is bounded by 
$$
\sum_{\substack{d\mid q\\d>s}} \sum_{c|d} 
 \frac{ c^{-1/2+\epsilon} }{ (d/c)^{1-\epsilon} } 
= 
\sum_{\substack{d\mid q\\d>s}} \frac {1}{d^{1-\epsilon}}
\sum_{c|d} c^{1/2+\epsilon}
$$
Now the sum $\sum_{c|d} c^{1/2+\epsilon}$ is bounded by 
$\tau(d) d^{1/2+\epsilon}\ll d^{1/2+\epsilon'}$, so the above is
bounded by  
$$
\sum_{\substack{d\mid q\\d>s}} d^{-1/2+\epsilon} \ll s^{-1/2+\epsilon}
$$
by Lemma \ref{bounding divisor sums}. 
This bounds the contribution of $c$, $\sGG$ with $c\supp(\sGG)>s$. 

If $c\disc(\sGG)>s$ then $\frac{s}{\disc(\sGG)} \leq c$. This,
together with \eqref{e:lattice-supp-sum} implies that
$$
\frac 1s \sum_{c\mid q} c^{-1/2+\epsilon}
\sum_{\substack{\supp(\sGG)\mid \frac qc\\ c\disc(\sGG)>s \\c\supp(\sGG)<s 
}}
\supp(\sGG)^\epsilon \frac{s  }{\disc(\sGG)} 
$$
$$
\ll
\frac 1s \sum_{c\mid q} c^{1/2+\epsilon}
\sum_{\substack{g\mid \frac qc \\ cg < s}}
g^\epsilon 
\ll
s^{-1/2+\epsilon} 
\sum_{\substack{c\mid q \\ c < s }} 
\sum_{\substack{g\mid q \\ g < s}}
1
$$
$$
\ll
s^{-1/2+\epsilon}
\left( \sum_{\substack{c\mid q \\ c < s }} 1 \right)^2
\ll
s^{-1/2+\epsilon},
$$ 
since $\sum_{\substack{c\mid q \\ c < s }} 1 \ll s^\epsilon$ by lemma
\ref{bounding divisor sums}. Consequently \eqref{truncated sum} is
$O(s^{-1/2+\epsilon})$. (Note that we only used $\supp(\sGG) \ll
s^{r(r-1)/2}$ to bound \eqref{extra term}!)


\subsection{The case $c\disc(\sGG)\leq s$} 
Fix $c\geq 1$ and $\sGG$ and partition the lattice points in 
$s\CC\cap L(\sGG)$ into two subsets as follows: Fix a {\em reduced}
fundamental cell (see \ref{reduced condition}) $P = P(\sGG)$ for the lattice
$L=L(\sGG)$. Then $cP$ is a reduced fundamental cell for the dilated
lattice $cL$. We can tile $\R^{r-1}$ by the translates $h_c+cP$, $h_c\in
cL$. 
\begin{defn}
We say that $x\in L\cap s\CC$ is {\em $c$-interior} if there is some
$y\in cL$  so that $x\in y+cP\subseteq s\CC$. We say that
$x\in L\cap s\CC$ is a {\em $c$-boundary} point otherwise.
\end{defn}
Note that the notion depends on $c$ and on the choice of a fundamental
cell $P$ for $L$. 

An important fact is that if $\dist(x, \partial(s\CC))\gg_r c\disc(L)$
then $x$ is $c$-interior. This follows from Lemma
\ref{diameter-discriminant bound} since
$\diam(cP) \ll_r c\disc(L)$.    

\begin{lem}\label{reps for Z/cZ} 
Let $P$ be a fundamental cell for the lattice $L\subseteq \Z^{r-1}$,
$c\geq 1$ so that  $\gcd(c,\disc(L)) =1$. Then for $y\in cL$, the
intersection $L\cap (y+cP)$ with $L$ of the translate of the dilated
cell $y+cP$   consists of a full set of representatives of
$\Z^{r-1}/c\Z^{r-1}$. 
\end{lem} 
\begin{proof}
If $P=\{\sum_{j=1}^{r-1} x_j \vl_j: 0\leq x_j<1 \}$ then the $c^{r-1}$
lattice points $y+\sum_{j=1}^{r-1} n_j \vl_j$, $n_j=0,1,\dots ,c-1$ in
$L\cap y+cP$ are clearly inequivalent modulo $cL$, and are the only
points of $L$ in this intersection. We will show that if
$\gcd(c,\disc(L)) =1$ then they are inequivalent modulo $c\Z^{r-1}$. 
To see this, it suffices to show that $L\cap c\Z^{r-1} = cL$.  
By the theorem on elementary divisors, there is a basis 
$\{\ve_j\}$ of $\Z^{r-1}$ and integers $d_j\geq 1$ so that
$\{d_j\ve_j\}$ is a basis of $L$, and $\disc(L) = \prod_{j=1}^{r-1}
d_j$. 
If $x\in L\cap c\Z^{r-1}$ then $x=\sum_{j=1}^{r-1} m_j d_j\ve_j \in L$
and also $x=c\sum_{j=1}^{r-1}n_j \ve_j \in c\Z^{r-1}$. Comparing
coefficients we find 
\begin{equation}\label{coefficients} 
 m_j d_j= c n_j, \qquad j=1,\dots, r-1
\end{equation}
Now since $d_j \mid \disc(L)$ and $\gcd(c,\disc(L)) = 1$, we have that
$\gcd(c,d_j) =1$ and so \eqref{coefficients} shows that $m_j=0\mod c$
and $x\in cL$.
\end{proof}

\begin{lem}\label{Lemma: count}
a) The number of points $y$ of $cL$ so that $y+cP\subset s\CC$ is 
$$
\frac{\vol(s\CC)}{\disc(cL)} + O\left((\frac sc)^{r-2}\right)
$$

b) The number of $c$-boundary points of $L$ is $\ll cs^{r-2}$.
\end{lem}
\begin{proof} 
If $y = cz\in cL$ then $y+cP\subseteq  s\CC$ if and only if 
$z\in L\cap\frac sc \CC$ and $z+P\subseteq \frac sc \CC$. 
Thus we need to count $N: = \#\{z\in L\cap\frac sc \CC: z+P\subseteq
\frac sc \CC\}$.  
An upper bound is obtained by a packing argument - since the
translates $z+P$ are disjoint and contained in $\frac sc\CC$, we get 
$$
N\vol(P) \leq \vol(\frac sc \CC)
$$
and so 
\begin{equation}\label{upper bound for N}
N\leq \frac{\vol(\frac sc \CC)}{\disc(L)} =
\frac{s^{r-1}\vol(\CC)}{c^{r-1}\disc(L)} 
\end{equation}
For a lower bound, note that if $z\in L\cap \frac sc\CC$ satisfies
$\dist(z,\partial (\frac sc\CC))>\diam(P) $ then $z+P\subseteq \frac
sc \CC$. By the Lipschitz principle (Lemma~\ref{Lipschitz lemma}) and
Lemma~\ref{convexity}, the number 
$\tilde N$ of such points  is 
$$
\tilde N = \frac{ \vol\{ x\in \frac sc \CC: \dist(x,\partial(\frac sc
\CC))\geq  \diam(P)\}} {\disc(L)} + 
O\left((\frac sc)^{r-2}\right) 
$$
Further, 
$$
\vol\{ x\in \frac sc \CC: \dist(x,\partial(\frac sc \CC))\geq \diam(P)\}
 =  \vol (\frac sc \CC) +
O\left(\diam(P)(\frac sc)^{r-2}\right)
$$
and so
$$
\tilde N = \frac{\vol (\frac sc \CC)}{\disc(L)} +
O\left(\frac{\diam(P)(\frac sc)^{r-2}}{\disc(L)} + (\frac sc)^{r-2}\right)
 = \frac{\vol (\frac sc \CC)}{\disc(L)} + O\left((\frac sc)^{r-2}\right)
$$
because $\diam(P) \ll_r \disc(L)$. 

Since $N\geq \tilde N$, together with the 
upper bound \eqref{upper bound for N} we find 
$$
N= \frac{\vol (\frac sc \CC)}{\disc(L)} + O\left((\frac sc)^{r-2}\right)
$$

b) For the number of $c$-boundary points, we subtract the number of 
$c$-interior points from the total number of points of $L\cap s\CC$. 
The total number of points in $L\cap s\CC$ is by the Lipschitz
principle (Lemma \ref{Lipschitz lemma}) 
\begin{equation}\label{total number}
L\cap s\CC = \frac{\vol (\ s \CC)}{\disc(L)} + O(s^{r-2}) 
\end{equation}
To count the number of $c$-interior points, we can write each uniquely
as $y+p$, with $y$ as in part a and $p\in L\cap cP$. 
Now $\#(L\cap cP) = c^{r-1}$ (see Lemma \ref{reps for Z/cZ}) and so by part a, 
the number $c$-interior points is 
\begin{equation}\label{number c-interior}
N c^{r-1} = \frac{\vol (\ s \CC)}{\disc(L)} + O(cs^{r-2})
\end{equation}
Subtracting \eqref{number c-interior} from \eqref{total number} gives part
b. 
\end{proof}

Fix $\sGG$, $c\geq1$ with $c\disc(\sGG)\leq s$. 
Note that since  $q$ is square-free and  $\supp(\sGG)\mid \frac qc$,
we have $\gcd(c,\disc(\sGG)) =1$.  
We now estimate the
sum 
$$
\sum_{h\in L(\sGG) \cap s\CC} a(h,c) \Delta(h,c)
$$
We divide this into two sums, $\Sigma_{int}$ over the $c$-interior
points and $\Sigma_{bd}$ over the $c$-boundary points. 
We use $a(h,c)\Delta(h,c) \ll c^{1/2+\epsilon}$ to bound $\Sigma_{bd}$
by: 
$$
 \#\{c\mbox{-boundary points}\}c^{1/2+\epsilon} 
\ll cs^{r-2} c^{1/2+\epsilon} = c^{3/2+\epsilon} s^{r-2}
$$
The contribution of the $c$-interior points is computed by writing
each 
such $h$ as $h=y+h_0$ with $h_0\in cP\cap L$ and $y\in cL \cap s\CC$. 
For each $y$ we get all possible $h_0$, which run over a full set of
representatives of $\Z^{r-1}/c\Z^{r-1}$ since $\gcd(c,\disc(\sGG)) =1$
(Lemma \ref{reps for Z/cZ}). Denote the
number of such $y$ by $N$; by Lemma~\ref{Lemma: count} part a, 
$N =  \frac{\vol (\frac sc \CC)}{\disc(L)} + O\left((\frac sc)^{r-2}\right)$.  
Moreover 
$$
a( y+h_0,c)\Delta(y+h_0,c) = a(h_0,c)\Delta(h_0,c)
$$
since $y\in cL(\sGG)\subset c \Z^{r-1}$. 
Thus 
\begin{equation*}
\begin{split}
\Sigma_{int} &= N \sum_{h_0 \mod c} a(h_0,c)\Delta(h_0,c)\\
& = \left( \frac{\vol (\frac sc \CC)}{\disc(L)} + O\left((\frac
sc)^{r-2}\right)\right) \sum_{h_0 \mod c} a(h_0,c)\Delta(h_0,c) \\
& = \frac{\vol (\frac sc \CC)}{\disc(L)}
\sum_{h_0 \mod c} a(h_0,c)\Delta(h_0,c) + 
O\left( (\frac sc)^{r-2} c^{r-1} c^{1/2+\epsilon}\right) \\
&= \frac{\vol (\frac sc \CC)}{\disc(L)}
\sum_{h_0 \mod c} a(h_0,c)\Delta(h_0,c) +  
O(c^{3/2+\epsilon} s^{r-2})
\end{split}
\end{equation*}
Thus the total contribution of the pairs with $c\disc(\sGG) \leq s$ is
\begin{multline}\label{cgamma<1}
\frac{s}{2^{r\omega(q)}} \sum_{c\mid q} \frac 1c
\sum_{\substack{\supp(\sGG)\mid \frac qc \\c\disc(\sGG) \leq s } }
\lambda(\sGG) 
\sum_{h\in s\CC\cap L(\sGG)}a(h,c) \Delta(h,c) \\
 = \frac{s}{2^{r\omega(q)}} \sum_{c\mid q} \frac 1c
\sum_{\substack{\supp(\sGG)\mid \frac qc \\c\disc(\sGG) \leq s }}
\lambda(\sGG) 
 \frac{\vol (s \CC)}{c^{r-1}\disc(\sGG)}
\sum_{h_0 \mod c} a(h_0,c)\Delta(h_0,c) \\
+O\left(  \frac{s}{2^{r\omega(q)}} \sum_{c\mid q} \frac 1c
\sum_{\substack{\supp(\sGG)\mid \frac qc \\c\disc(\sGG) \leq s }}
|\lambda(\sGG)|
c^{3/2+\epsilon} s^{r-2} \right)
\end{multline}

To estimate the error in \eqref{cgamma<1}, note that the condition
$c\disc(\sGG) \leq s$ 
implies $c\supp(\sGG) \leq s$ since $\supp(\sGG) \leq \disc(\sGG)$, so
for an upper bound we may replace the summation over pairs satisfying
the former condition by the sum over pairs satisfying the latter; this
gives 
%
%
(noting that $2^{\omega(q)} \geq s$)
\begin{equation*}
\begin{split}
\frac{s}{2^{r\omega(q)}} \sum_{c\mid q} \frac 1c
\sum_{\substack{\supp(\sGG)\mid \frac qc \\c\disc(\sGG) \leq s } }
|\lambda(\sGG)| c^{3/2+\epsilon} s^{r-2} 
&\ll
s^{-1+\epsilon} \sum_{c\mid q} c^{1/2+\epsilon} 
\sum_{\substack{\supp(\sGG)\mid \frac qc \\c\supp(\sGG) \leq s }  }
|\lambda(\sGG)| \\
&\ll s^{-1+\epsilon} \sum_{c\mid q} c^{1/2+\epsilon} 
\sum_{\substack{g \mid \frac qc \\cg \leq s }}
\sum_{\supp(\sGG) = g} |\lambda(\sGG)|
\end{split}
\end{equation*}
Now $|\lambda(\sGG)| \ll\supp(\sGG)^\epsilon$ and the number of $\sGG$
with 
$\supp(\sGG) = g$ is $O(g^\epsilon)$, which is
$O(s^\epsilon)$ since $g\leq cg\leq s$, so that the above is
bounded by 
\begin{equation*}
s^{-1+\epsilon} \sum_{c\mid q} c^{1/2+\epsilon} 
\sum_{\substack{g \mid \frac qc \\cg \leq s }}1
\end{equation*}
The number of small divisors $g$ of $\frac qc$ with $g
\leq s/c\leq s$ is at most $s^\epsilon$, so the above is at most 
$$
s^{-1+\epsilon} \sum_{\substack{c\mid q\\c\leq s}} c^{1/2+\epsilon} 
\ll s^{-1+\epsilon}s^{1/2+\epsilon} \#\{c\mid q: c\leq s\} \ll
s^{-1/2+\epsilon'} 
$$
which gives that the error term in \eqref{cgamma<1} is $O(s^{-1/2+\epsilon})$.

We now extend the sum of the first term in \eqref{cgamma<1} to all the
pairs $c$, $\sGG$, introducing an error which was bounded in 
section \ref{sec cgamma>s} by $O(s^{-1/2+\epsilon})$. 
(This is the term \eqref{truncated sum} which was bounded without using
the condition $\supp (\sGG)\ll s^{r(r-1)/2}$.)

In summary we find that 
\begin{prop}
\begin{multline}\label{R=M+E}
R_r(\CC,q) = 
  \frac{s}{2^{r\omega(q)}} \sum_{c\mid q} \frac 1c
\sum_{\supp(\sGG)\mid \frac qc }
\lambda(\sGG) 
 \frac{\vol (s \CC)}{c^{r-1}\disc(L)}
\sum_{h_0 \mod c} a(h_0,c)\Delta(h_0,c) \\
+O(s^{-1/2+\epsilon})
\end{multline}
\end{prop}

\newpage
\section{The main term}\label{sec:mainterm}
We now treat the main term of \eqref{R=M+E}. Define 
$$ \MM = 
\frac{s}{2^{r\omega(q)}} \sum_{c\mid q} \frac 1c
\sum_{\supp(\sGG)\mid \frac qc }
\lambda(\sGG) 
 \frac{\vol (s \CC)}{c^{r-1}\disc(\sGG)}
\sum_{h_0 \mod c} a(h_0,c)\Delta(h_0,c)
$$ 
we will show 
$$
\MM= \vol(\CC)
$$ 
which with  \eqref{R=M+E} will prove Theorem \ref{rlevel thm}. 

The sum over $h \mod c$ is multiplicative: 
$$
\sum_{h\mod c} a(h,c)\Delta(h,c)=  \prod_{p\mid c} \sum_{h\mod p}
a(h,p)\Delta(h,p)
$$
Furthermore, by Lemma \ref{sum a(h,p)Delta(h,p)}
\begin{equation*}
\sum_{h\mod p} a(h,p)\Delta(h,p) 
= (p+1)^r - p^r\sum_{\sGp} \lambda(\sGp)
p^{-\codim H_{\sGp}} 
\end{equation*}
Now note that 
since $p^{\codim H_{\sGp}}  = \disc(\sGp)$ we get 
\begin{equation*}
\begin{split}
\MM & =
\frac{s}{2^{r\omega(q)}} \sum_{c\mid q} \frac 1c
\sum_{\supp(\sGG)\mid \frac qc} \lambda(\sGG)  
\frac{s^{r-1}\vol(\CC)}{c^{r-1}\disc(\sGG)} 
\prod_{p\mid c}\left(
(p+1)^r - p^r\sum_{\sGp}\frac{\lambda(\sGp) }{\disc(\sGp)}
 \right) \\
& = \frac{\vol(\CC)s^r}{2^{r\omega(q)}} \sum_{c\mid q} \frac 1{c^r} 
\sum_{\supp(\sGG)\mid \frac qc}\frac{\lambda(\sGG)}{\disc(\sGG)}  
\prod_{p\mid c}\left( 
(p+1)^r - p^r\sum_{\sGp}\frac{\lambda(\sGp) }{\disc(\sGp)}
\right) 
\end{split}
\end{equation*}
Further, 
$$
\sum_{\supp(\sGG)\mid \frac qc} \frac{\lambda(\sGG)}{\disc(\sGG)}  = 
\prod_{p\mid \frac qc} \sum_{\sGp} \frac{\lambda(\sGp)}{\disc(\sGp)}
$$
Therefore we find that 
\begin{equation*}
\begin{split}
\MM &= \vol(\CC)\frac 1{\sigma_{-1}(q)^r} \sum_{c\mid q} 
\prod_{p\mid \frac qc} \sum_{\sGp} \frac{\lambda(\sGp)}{\disc(\sGp)} 
\prod_{p\mid c}\left( 
(1+\frac 1p)^r - \sum_{\sGp}\frac{\lambda(\sGp) }{\disc(\sGp)} \right) \\
&= \vol(\CC)\frac 1{\sigma_{-1}(q)^r} \sum_{c\mid q}
A(\frac qc) B(c)
\end{split}
\end{equation*}
Thus $\MM$ is a multiple of the Dirichlet convolution of the
multiplicative functions  $A$, $B$, with $A(1) = B(1) = 1$, 
$$
A(p)  = \sum_{\sGp} \frac{\lambda(\sGp)}{\disc(\sGp)} 
$$
and (since $(1+1/p)^r = \sigma_{-1}(p)^r$)  
\begin{equation}\label{A and B}
B(p)  = \sigma_{-1}(p)^r - A(p) 
\end{equation}
Now we have 
$$
(A\ast B)(q):= 
\sum_{c\mid q}A(\frac qc) B(c)  = \prod_{p|q} \left( A(1)B(p) + A(p)B(1)
\right)  =\prod_{p\mid q} \sigma_{-1}(p)^r = \sigma_{-1}(q)^r
$$ 
by \eqref{A and B}. This finally gives the main term of $R_r(\CC,q)$: 
$$
\MM = \vol(\CC) \frac 1{\sigma_{-1}(q)^r}(A\ast B)(q) = 
 \vol(\CC) \frac 1{\sigma_{-1}(q)^r}\sigma_{-1}(q)^r = \vol(\CC)
$$ 

\newpage
\appendix
\section{Recovering the level spacing from the
correlations}\label{app:spacings}  

In this appendix, we explain how to recover the various spacing
distributions from the correlation functions. This is
well-known in the physics literature (e.g. \cite{Mehta})   and is
certainly implicit in Hooley's work \cite{Hooley2, Hooley3,
Hooley-survey}, but we do not know of a good source for it in the
mathematical literature. A very detailed treatment of this and more will
appear in a forthcoming book by Katz and Sarnak \cite{KS}. 

We begin with $\R/\Z$ 
which we think of as the circle with  unit circumference. 
We denote by $\{x\}$  the fractional part of $x$: If $n\leq x <n+1$, $n$
integer,  then $\{x\}=x-n$.   We set 
$$
((x)) = \begin{cases} \{x\}& 0\leq \{x\}< \frac 12 \\
                      \{x\}-1&\frac 12\leq \{x\}<1 \end{cases}
$$ 
We will order the points in $\R/\Z$  counter-clockwise and write
$x\succ y$ if the points lie in a segment of length $<1/2$ on $\R/\Z$ 
and $x$ follows  $y$.  
The {\em signed distance} on $\R/\Z$ is given by $((x-y))$; 
thus $-1/2\leq ((x-y)) < 1/2$. 
In terms of the signed distance, $x\succ y$ if and only if $((x-y))>0$.

Given a finite set $S$ of $N$ points on $\R/\Z$, and 
$k\geq 2$, the {\em $k$-level correlation functions} measure clustering
properties of the sequence $S\subset \R/\Z$ on a scale of the mean 
spacing $1/N$: For a $k$-tuple of points $x=(x_1,\dots,x_k)$ of
$S$, the oriented distance vector is 
\begin{equation}\label{orientedD}
D(x) = \left( ((x_1-x_2)),\dots,((x_{k-1}-x_k)) \right)
\end{equation}
Given a bounded set $\CC \subset \R^{k-1}$, we define the $k$-level
correlation as 
$$
R_k(\CC,S) = \frac 1N \#\{ x\in S^k: D(x) \in \frac 1N \CC \}
$$

As an example, let $\Delta^{k-1} \subset \R^{k-1}$ be 
the standard open simplex  
$$
\Delta^{k-1} = 
\{ (y_1,\ldots y_{k-1}) \ | \ y_i >0, \sum_{i=1}^{k-1} y_i < 1 \}
$$ 
and for $t>0$ set $\CC=t \Delta^{k-1}$. Then if $N>2t$, $D(x)\in
\frac 1N\CC = \frac tN\Delta^{k-1}$ means that 
\begin{enumerate}
\item $((x_i-x_{i+1}))>0$, that is $x_1\succ x_2\succ\dots \succ x_k$;
\item The points all lie in an arc of length at most $t/N$. 
\end{enumerate} 

As another example, write $k-1=i+j$ and for $t_1, t_2>0$ set 
$\CC= t_1\Delta^i\times t_2 \Delta^j $, which we can write as  
$$
\CC =
\{(y_1,\dots y_k):y_m>0, y_1+y_2+\dots +y_i< t_1, y_{i+1}+\dots +
y_{i+j}< t_2 \} 
$$  
Then $D(x) \in \frac 1N \CC$ iff $x_1\succ x_2\succ\dots \succ x_k$
and $x_1,\dots, x_{i+1}$ lie in an arc of length $<t_1/N$, and
$x_{i+1},\dots, x_{i+j+1}=x_k$ lie in an arc of length $<t_2/N$.  

Given any subset $T\subseteq S$ which is contained
in a semi-circle, 
the ordering gives us unique initial and final elements of $T$, and we
can write $T=\{x_{init}=x_1\prec x_2\prec \dots\prec x_{fin} \}$. 
We denote by $|T|$ 
the number of elements of $T$, and by $\diam(T)$ the distance
$\dist(x_{init},x_{fin})$ between the initial and final points of
$T$. 
If $T$ consists of just the initial and final points, we say that
$T$ is a {\em consecutive pair}. A {\em consecutive $k$-tuple} of $S$ 
is a $k$-tuple of elements 
$x_1=x_{init}\prec \dots \prec x_k = x_{fin}$  so that there are no points of
$S$ between $x_j$ and $x_{j+1}$, for $1\leq j<k$.

For $x<1/2$, let $N_k(x)$ be the number of  $k$-tuples of
diameter smaller than $x$; this is zero if $k\gg1$. 
It is clear from the definitions and the discussion above that we can
describe these functions in terms of the correlation function of the
simplex $x \Delta^{k-1}$ by
\begin{equation}\label{N intermsof R}
R_k(x\Delta^{k-1},S) = \frac 1N N_k(\frac xN)
\end{equation}
Further, let $g(x)$ be the number of {\em consecutive pairs} of diameter
less than $x$, that is the number of spacings between consecutive
elements of $S$ of length less than $x$. 
We may express $g$ in terms of an alternating sum of
$N_k$'s as follows:  
\begin{lem}
  \label{l:sandwich}
  With $g$ and $N_k$ as above, we have for $x<1/2$ 
  $$
  g(x) = \sum_{k \geq 2} (-1)^{k} N_k(x).
  $$ 

  Moreover, for all $n \geq 1$, we have the inequalities
  $$
  \sum_{k=2}^{2n+1} (-1)^{k} N_k(x)
  \leq   g(x) \leq
  \sum_{k=2}^{2n} (-1)^{k} N_k(x)
  $$
\end{lem}

Before giving the proof, 
we will need the following elementary lemma on sums of binomial coefficients. 
\begin{lem}
\label{l:alternating-binomial}
Let $m \geq 0$ be an integer. Then $\sum_{i=0}^m (-1)^i \binom{m}{i} =
0$ unless $m=0$, in which case the sum equals $1$. 
Moreover, 
$$
\sum_{i=0}^{2n+1} (-1)^i \binom{m}{i} \leq 
\sum_{i=0}^m (-1)^i \binom{m}{i} \leq
\sum_{i=0}^{2n} (-1)^i \binom{m}{i}.
$$
\end{lem}
\begin{proof} 
The first part is just the binomial expansion of $(1-1)^m$. As for the
second part,  if $m\geq 1$ use  
the identity $\binom{m}{i} =\binom{m-1}{i}+ \binom{m-1}{i-1}$ to find 
$\sum_{i=0}^k (-1)^i \binom{m}{i} = (-1)^k \binom{m-1}{k}$ from which
the claim follows. 
\end{proof}
We can now prove lemma~\ref{l:sandwich}:

\begin{proof} 
 For each  pair $T=\{a\succ b\}$, of diameter less than
$1/2$, we associate $X_T$, the set of all  $i$-tuples 
$x_1\succ \ldots\succ x_i$ in $S$ such that $(x_1,x_i)=(a,b)$. The set of all 
tuples of diameter less than $x$ is thus expressed as a {\em disjoint
  union} of the $X_T$'s as $T$ ranges over all  pairs of
diameter less than $x$. If we let $N_i^T$ be the number of $i$-tuples in
$X_T$ then  $N_i = \sum_T N_i^T$. But $N_i^T = \binom{|T|-2}{i}$,
so by lemma~\ref{l:alternating-binomial}, 
$
\sum_{i \geq 2} (-1)^i N_i^T
$
is zero unless $T$ is a {\em consecutive} pair, in which case the
alternating sum  is one. Summing over all consecutive pairs we get
that  $g(x) = \sum_{k \geq 2} (-1)^{k} N_k(x).$
Lemma~\ref{l:alternating-binomial} also gives that for $n>0$, 
$$
\sum_{i=2}^{2n+1} (-1)^i N_i^T
\leq \sum_{i \geq 2} (-1)^i N_i^T \leq
\sum_{i=2}^{2n} (-1)^i N_i^T 
$$
Summing over all $T$ we get the second assertion. 
\end{proof}

\subsection{The joint level spacing}
An  $(i,j)$-tuple of diameter $(x,y)$ is  
an  $(i+j)$-tuple $x_1\succ \ldots\succ x_i \succ
x_{i+1}\succ \ldots\succ x_{i+j}$ (all lying in an arc of length $<1/2$)  
such that $\dist(x_i,x_1) =x$ and $\dist(x_{i+j},x_i) =y$. 

For $i \geq 2$,  $j \geq 1$ and $x+y<1/2$ 
we let $N_{i,j}(x,y)$ be the number of
$(i,j)$-tuples of diameter at most $(x,y)$. Let $g(x,y)$ be
the number of consecutive triples $x_1\succ x_2\succ x_3$ of diameter smaller
than $(x,y)$. Analogously to lemma~\ref{l:sandwich} we have:  
\begin{lem}
If we let $A_k(x,y) = \sum_{i+j=k} N_{i,j}(x,y)$, then
$$
g(x,y) = \sum_{k \geq 3} (-1)^{k+1} A_k(x,y).
$$
Moreover, for $n \geq 0$ we have the inequalities
$$
\sum_{k=3}^{3 + 2n +1} (-1)^{k+1} A_k(x,y)
\leq g(x,y) \leq 
\sum_{k=3}^{3 + 2n} (-1)^{k+1} A_k(x,y).
$$
\end{lem}
\begin{proof} For each triple $T=\{a\succ b\succ c\}$ of diameter at
most $(x,y)$, 
let $X_T$ be the set of $(i,j)$-tuples
$x_1\succ  \ldots\succ  x_i\succ  x_{i+1}\succ \ldots \succ x_{i+j}$ 
such that $(x_1,x_i,x_{i+j})=(a,b,c)$, and let $N^T_{i,j}$ be the
number of 
$(i,j)$-tuples in $X_T$. We may write the set of $(i,j)$-tuples of
diameter smaller than $(x,y)$ as a {\em disjoint} union of $X_T$'s, as
$T$ ranges over all $(2,1)$-tuples with diameter at most
$(x,y)$. Given $T$, we may count tuples of type $(i,j)$ in 
$X_T$ as follows: Let $M,N$ be the number of elements of $S$ 
{\em between} $a,b$ and $b,c$ respectively (we allow both $M$ and $N$
to be zero.) Then $N^T_{i,j} = \binom{M}{j-2} \binom{N}{i-1}$. Moreover, 
$A^T_k =\sum_{i+j=k} N^T_{i,j} = \binom{M+N}{k-3}$ since there are
$\binom{M+N}{k-3}$ ways of choosing $k-3$ objects out of $M$ ``blue'' and
$N$ ``red'' objects. By lemma ~\ref{l:alternating-binomial}, we see that 
$
\sum_{k \geq 3 } (-1)^{k+1} A^T_k =
\sum_{k \geq 3} (-1)^{k+1} \binom{M+N}{k-3}
$
is zero unless $T$ is a consecutive $(2,1)$-tuple, in which case it is
one. Now, lemma~\ref{l:alternating-binomial} together with $A^T_k =
\binom{M+N}{k-3}$ shows that 
$$
\sum_{k=3}^{3+2n-1} (-1)^{k+1} A^T_k
\leq \sum_{k \geq 3} (-1)^{k+1} A^T_k \leq
\sum_{k=3}^{3+2n} (-1)^{k+1} A^T_k.
$$
Summing over all triples $T$ of diameter at most $(x,y)$
we are done. \end{proof}

\subsection{Applications to squares mod $q$.}

We let 
$$
S_q=\{ \frac nq: 0\leq n\leq q-1, \quad n \mbox{ a square modulo }q
\} \subset \R/\Z
$$
be the image in $\R/\Z$ of the set of squares in $\Z /q \Z$.  
The mean spacing between elements of $S_q$ is $1/N_q$, where $N_q$ is
the number of squares modulo $q$. For $x>0$, $g_q(\frac x{N_q})$
is the number of consecutive pairs in $S_q$ of diameter  at most
$x/N_q$, that is the number of normalized consecutive spacings of
length $<x$, and we set 
$$
\tilde P(x) = \lim_{q\to \infty} \frac 1{N_q} g_q(\frac x{N_q})
$$
This is the limiting proportion of normalized consecutive spacings in
$S_q$ of length at most $x$ (this normalization sets the mean spacing
to be unity). $\tilde P(x)$ is the cumulant of the level spacing
distribution $P(s)$ of the introduction. 
Likewise we set for $x,y>0$, 
$$
\tilde P(x,y) =  \lim_{q\to \infty} 
\frac 1{N_q} g_q(\frac x{N_q},\frac y{N_q} ) 
$$
the cumulant of the joint level spacing distribution. 



For a bounded convex set $\CC\subset \R^{k-1}$, not intersecting the
walls, and $N\gg1$, $\frac 1{N_q}\CC$ will be contained in the cube
$(-1/2,1/2)^{k-1}$. For $x=\frac nq \in S_q^k$, ($0\leq n_i<q$ are
squares modulo $q$) the oriented
distance vector $D(x)$ (see \eqref{orientedD}) will lie in $\frac 1{N_q} \CC$
if and only if 
there is an  integer vector $h\in \frac q{N_q}\CC \cap \Z^{k-1}$ so
that 
$$
x_i-x_{i+1} = h_i \mmod q , \quad 1\leq i \leq k-1
$$
Denoting by $N(h,q)$ the number of solutions of the above system in
squares $n_i$ modulo $q$, we have found that the correlation function 
$R_k(\CC,q) := R_k(\CC,S_q)$ satisfies 
\begin{equation}
R_k(\CC,q)  = \frac 1{N_q} \sum_{h\in s\CC\cap \Z^{k-1}} N(h,q)
\end{equation}
with $s=q/N_q$. 

\begin{lem}
If $x,y>0$ then
$$
\tilde P(x) = 1 - e^{-x}
$$
and 
$$
\tilde P(x,y) = (1 - e^{-x})(1-e^{-y})
$$
\end{lem}
\begin{proof} 

As noted above (see \eqref{N intermsof R}), we can express the functions 
$N_k(x)$ in terms of the correlation functions associated to the
simplex $x\Delta^{k-1}$, whose volume is $\frac {x^{k-1}}{(k-1)!}$: 
$$
R_k(x\Delta^{k-1};q) = \frac 1{N_q} N_k(\frac x{N_q})
$$
 From theorem~\ref{rlevel thm} we know that 
\begin{equation*}
\begin{split}
R_k(x \Delta^{k-1};q) &=
x^{k-1} \vol(\Delta^{k-1}) +O_k(s^{-1/2+\epsilon}) \\
& = 
\frac{x^{k-1}}{(k-1)!} +O_k(s^{-1/2+\epsilon})
\end{split}
\end{equation*}
By lemma~\ref{l:sandwich}
we see that for $n>0$,
$$
\sum_{i=1}^{1+2n+1}(-1)^{i+1} \frac{x^i}{i!} \leq 
\liminf_{q \rightarrow \infty} \frac{g_q(\frac x{N_q})}{N_q} 
$$
and
$$
\limsup_{q \rightarrow \infty}\frac{g_q(\frac x{N_q})}{N_q}  
\leq \sum_{i=1}^{1+2n}(-1)^{i+1} \frac{x^i}{i!}
$$
Letting $n \rightarrow \infty$ and noting that the above polynomials
are truncations of the Taylor series of $1-e^{-x}$ we are done.

For the second part of the lemma, recall that 
$N_{i,j}(x,y)$ is the number of ordered $i+j$-tuples of elements of
$S_q$ such that the first $i$ are contained in an interval of length
$x$, and the last $j$ elements lie in an interval of length $y$.  
Thus, analogously to \eqref{N intermsof R}, $N_{i,j}(x,y)$ is a scaled
version of the 
$(i+j-1)$-correlation with 
respect to the convex set $x \Delta^{i-1} \times y \Delta^{j}$:  
$$
\frac{N_{i,j}(\frac x{N_q},\frac y{N_q})}{N_q}  = 
R_{i+j}( x \Delta^{i-1} \times y \Delta^{j} ; q).
$$
By Theorem \ref{rlevel thm}, 
$$
R_{i+j}( x \Delta^{i-1} \times y \Delta^{j} ; q)
=
\frac{x^{i-1} y^j}{(i-1)! j!} + O_{i,j}(s^{-1/2+\epsilon})
$$
since
$$
\vol( x \Delta^{i-1} \times y \Delta^{j}) 
= \frac{x^{i-1} y^j}{(i-1)! j!} 
$$
Letting $A_k(x,y)=
\sum_{i+j=k} N_{i,j}(x,y)$ and using
lemma~\ref{l:alternating-binomial}, we get
$$
\limsup_{q \rightarrow \infty} \frac{1}{N_q}g_q(\frac x{N_q},\frac y{N_q}) 
\leq
\sum_{k=3}^{3+2n} (-1)^{k+1} 
\sum_{\substack{ i+j=k \\ i>1 \\ j>0 } }
\frac{x^{i-1} y^j}{(i-1)! j!}
$$
and
$$
\sum_{k=3}^{3+2n} (-1)^{k+1} 
\sum_{\substack{ i+j=k \\ i>1 \\ j>0 } }
\frac{x^{i-1} y^j}{(i-1)! j!}
\leq
\liminf_{q \rightarrow \infty} \frac{1}{N_q}g_q(\frac x{N_q},\frac y{N_q})
$$
Since the above polynomials are truncations of the 
Taylor series for $(1-e^{-x})(1-e^{-y})$, we are done. \end{proof}

\newpage
\section{Some Geometry of Numbers}\label{sec:geom}

\subsection{} 
Given a basis $\vl_1\dots,\vl_n$ of a lattice $L$ in $\R^n$, 
the {\em fundamental cell} is the half-open set 
$$
P(\{\vl_i\}):= \{ x_1\vl_1+\dots + x_n \vl_n: 0\leq x_i <1 \}
$$
It serves as a fundamental domain for the action of $L$ on $\R^n$ by
translations. The volume of $P(\{\vl_i\})$ is the {\em discriminant}
$\disc(L)$ of the lattice $L$:
$$
\vol( P(\{\vl_i\})) = |\det(\vl_1,\dots,\vl_n)| = \disc(L)
$$

\subsection{} 
We need the following basic fact (due to Mahler and Weyl) from
reduction theory: In any dimension $n\geq 1$, there are constants
$0<c_n'<c_n''$ so that any lattice $L\subset \R^n$ has a basis
$\vl_1\dots,\vl_n$ which is {\em reduced} in the sense that 
\begin{equation}\label{reduced condition}
c_n' \leq \frac{|\vl_1|\dots |\vl_n|}{\disc(L)} \leq c_n''
\end{equation}
This is a consequence of Minkowski's second theorem on successive
minima; see (Cassels \cite{cassels-geom}, Lemma V.8, p. 135) or 
(Siegel \cite{siegel}, X \S 6). This basis is not unique in general. 

\subsection{} 
We define the {\em diameter} $\diam(L)$ of the lattice $L$ to be the
minimum of the diameters of all fundamental cells for $L$.
\begin{lem}\label{diameter-discriminant bound}
The diameter of an {\em integer} lattice $L \subseteq \Z^n$ is bounded
by the discriminant of $L$: 
\begin{equation}
\diam(L) \ll_n \disc(L)
\end{equation}
the implied constant depending only on the dimension $n$.  
\end{lem}
\begin{proof}
It suffices to show that if $P(\{\vl_i\})$ is the fundamental cell of
an {\em integer} lattice  
$L \subseteq \Z^n$ with respect to a {\em reduced} basis $\{\vl_i\}$, then
the {\em diameter} of $P(\{\vl_i\})$ is bounded by the discriminant of
$L$: 
\begin{equation}
\diam(P(\{\vl_i\})) \ll_n \disc(L)
\end{equation}

To see this, note that since $L\subseteq \Z^n$ is an {\em integer} lattice, 
the length of any non-zero vector in $L$ is at least $1$, and then
this implies that a reduced basis has {\em bounded eccentricity}: 
\begin{equation}
1\leq |\vl_1| \leq |\vl_2| \leq \dots \leq |\vl_n| \leq c_n'' \disc(L)
\end{equation}
(assuming we ordered the basis vectors according to their length). 
Indeed, using (\ref{reduced condition}) together with $|\vl_i|\geq 1$ we
get an upper bound  for the longest basis vector $\vl_n$
$$
|\vl_n|= 1\cdot |\vl_n| \leq |\vl_1|\cdot |\vl_2|\cdot \dots \cdot
|\vl_n| \leq c_n'' \disc(L)
$$
Thus the diameter of the fundamental cell $P(\{\vl_i\})$ is at most 
$$
\sum_{i=1}^n |\vl_i| \leq n|\vl_n| \leq c_n'' \disc(L)
$$
as required. 
\end{proof}

\subsection{} 
It will be useful to note that for integer {\em dilates} $cL$ of a
lattice $L$, $c\geq 1$, the diameter scales {\em linearly}: $\diam(cL)
= c\diam(L)$, while the discriminant scales with $c^n$: $\disc(cL) =
c^n\disc(L)$. Thus to bound the diameter of a dilate of an integer lattice
we use 
\begin{equation}
\diam(cL) \ll_n c\disc(L)
\end{equation}

\subsection{The Lipschitz principle} 
\begin{defn} A set $\CC\subset \R^n$ is {\em of class} $m$ if the
intersection of every line with $\CC$ consists of at most $m$ intervals
(including the degenerate case when some of the intervals are
points), and if the same is true for the projection of $\CC$ on every
linear subspace. 
\end{defn}
Thus for instance a {\em convex} set is of class 1.

We will use the following form of the ``Lipschitz principle'' from the
geometry of numbers to estimate the number of lattice points in a
region of $\R^n$: 
\begin{lem} \label{Lipschitz lemma}
Let $L\subset \Z^n$ be an integer lattice of discriminant $\disc(L)$, 
and $\CC \subset \R^n$ a set of class $m$ (e.g. a convex set). Suppose
that $\CC$ lies in a ball of radius $R$ around the origin. 
Then 
\begin{equation}
\# (L\cap \CC) = \frac{\vol(\CC)}{\disc(L)} +O(R^{n-1})
\end{equation}
\end{lem}
This follows from the Lipschitz principle for the integer lattice 
proven by Davenport \cite{Davenport}, as adapted by W. Schmidt
(\cite{Schmidt}, Lemma 1).  

We will apply the Lipschitz principle to certain subsets of convex sets. 
For this purpose we will need: 
\begin{lem} \label{convexity}
Let $\CC \subset \R^n$ be a convex set, $d>0$ and define 
$$
\CC_d:= \{x\in \CC: \dist(x,\partial \CC) \geq d \}
$$
to be the set of points of $\CC$ of distance at least $d$  from the
boundary $\partial \CC$ of $\CC$. Then $\CC_d$ is convex. 
\end{lem}
\begin{proof}
What we need to show is that for any 
$x_1,x_2 \in \CC_d$, and   $\lambda \in [0,1]$, 
the point $x_3 = x_1 + \lambda (x_2-x_1)$ also lies in $\CC_d$, 
that is if $|y|\leq d$ then 
$ x_3+y\in \CC$. But $x_3+y = (x_1+y) + \lambda ( (x_2+y) - (x_1+y) )$,
i.e. $x_3+y$ lies on a line between $x_1+y$ and $x_2+y$. These two points 
lie in $\CC$ since $x_1,x_2\in \CC_d$. 
By convexity so does $x_3+y$. 
\end{proof}

\newpage

\newpage
\section{Counting small divisors} \label{app:divisors}

In the paper, 
we  need to use some estimates  for the number of divisors of $q$
that are smaller than  a fixed power of the mean spacing $s$.   
As is well known, the number of all divisors of $q$ is $O(q^\epsilon)$
for all $\epsilon>0$. This is not enough for our purposes, as we need
a bound which is $O(s^\epsilon)$. This is provided by the following
lemmas: 
\begin{lem}
\label{small divisors}
Let $q$ be square-free, $s=2^{\omega(q)}/\sigma_{-1}(q)$. Fix
$\alpha>0$. Then as $s\to \infty$  
$$
\#\{ d\mid q : d<s^\alpha\} = O(s^{\epsilon}). 
$$
for all $\epsilon>0$.
\end{lem} 
\begin{proof} We start by bounding products of $k$ distinct primes below by
$k^k$; we may assume that the primes are the first $k$
primes. Then by the Prime Number Theorem, 
$$
\log \prod_{i=1}^k p_i = 
\sum_{i=1}^k \log p_i \sim  p_k \sim k\log k
$$
Exponentiating we see that the
product is bounded below by $k^k$. Now,
$$
\#\{ d\mid q : d<s^\alpha \} =
\sum_j a_j,
$$ 
where $a_j=a(j,s^\alpha,q)$ is the number of divisors 
of $q$ that are smaller than $s^\alpha$ and have precisely $j$ prime
factors. But if $j>N$, where $N$ is the
smallest integer such that $N^N \geq s^\alpha$, then
$a_j=0$. Moreover, setting $w=\omega(q)$, we see that $a_j \leq
\binom{w}{j}$. Hence 
$$
\sum_{ \substack{ d\mid q \\ d<s^\alpha  } } 1 \leq \sum_{j \leq N}
\binom{w}{j} \leq N \binom{w}{N}.
$$
By Stirling's formula, $ \binom{w}{N} \ll \frac{ w^N}{(N/e)^N}$. Thus 
\begin{equation*}
\sum_i a_j \leq
N \binom{w}{N} \ll
 N \left( \frac{ w e }{ N } \right)^N  \ll
 N \left( \frac{ N \log(N) e }{ \alpha N \log(2)} \right)^N  
\end{equation*}
since $N^N \geq s^\alpha \gg 2^{w \alpha (1-\epsilon)}$ implies
that $w \leq \frac{N \log(N)}{\alpha \log(2)}$. 
Thus 
\begin{equation*}
\{ d\mid q : d<s^\alpha \}
\ll  N \left( \frac{ \log(N) e }{ \alpha \log(2) } \right)^N  \ll (C\log N)^N
\end{equation*}
But the last 
term is clearly  $O(s^\epsilon)$.  
\end{proof}

\begin{lem}\label{bounding divisor sums}
If $\alpha>0$ then 
$\sum_{\substack{d\mid q\\ d>s}} d^{-\alpha} \ll s^{-\alpha+\epsilon}$ 
\end{lem}
\begin{proof}
We divide the sum into two parts: One over $s<d<s^R$ and the other over
$d>s^R$ ($R$ is a parameter chosen later). For the first, we use the
fact that there are few (namely $ O(s^\epsilon)$)  
divisors $d$ of $q$ with $d<s^R$ to bound that contribution by 
$$
\sum_{\substack{d\mid q\\ s<d<s^R}}d^{-\alpha} \ll 
\sum_{\substack{d\mid q\\s<d<s^R}} s^{-\alpha} \ll
s^{-\alpha+\epsilon} 
$$
For the summands with $d>s^R$, use $d^{-\alpha} < s^{-R\alpha}$ and
$\tau(q) = 2^{\omega(q)} \ll s^{1+\epsilon}$ 
to get 
$$
\sum_{\substack{d\mid q\\ d>s^R}} d^{-\alpha} \ll s^{-R\alpha} \tau(q) \ll
s^{1-R\alpha+\epsilon}  
$$
Now choose $R>0$ so that $1-R\alpha <-\alpha$ to conclude the lemma.
\end{proof}

\newpage
\bibliographystyle{amsplain}

\end{document}